\documentclass[11pt]{article}
\usepackage{titlesec}
\usepackage{amscd,amssymb,stmaryrd}
\usepackage[fleqn]{amsmath}
\usepackage{graphicx}
\usepackage{subcaption}
\usepackage{comment}
\usepackage[utf8]{inputenc}
\usepackage[export]{adjustbox}
\usepackage{wrapfig}
\usepackage{amstext}
\usepackage{pstricks,pst-node,pst-plot,pst-coil}
\usepackage{amsthm}
\usepackage{amsmath}
\usepackage{color}
\usepackage{mathtools}
\usepackage{hyperref}
\usepackage[english]{babel}
\usepackage{booktabs}
\usepackage{mathrsfs}
\usepackage{float} 

\newcommand\norm[1]{\left\|#1\right\|}
\newcommand\abs[1]{\lvert#1\rvert}

\newcommand{\tforall}{\text{ for all }}

\newcommand{\eps}{\varepsilon}
\theoremstyle{definition}

\theoremstyle{remark}

\DeclareMathOperator*{\argmin}{argmin}
\title{A multi-stage deep learning based algorithm for multiscale model reduction}

\author{Eric Chung\thanks{Department of Mathematics, The Chinese University of Hong Kong, Shatin, Hong Kong}\ , \quad Wing Tat Leung\thanks{Department of Mathematics, University of California, Irvine, CA 92697, USA}\ , \quad Sai-Mang Pun\thanks{Department of Mathematics, Texas A\&M University, College Station, TX 77843, USA}\ , \quad and \quad Zecheng Zhang\thanks{Department of Mathematics, Texas A\&M University, College Station, TX 77843, USA}}

\date{\today}

\begin{document}
\maketitle 

\begin{abstract}
In this work, we propose a multi-stage training strategy for the development of deep learning algorithms applied to problems with multiscale features. Each stage of the proposed strategy shares an (almost) identical network structure and predicts the same reduced order model of the multiscale problem. The output of the previous stage will be combined with an intermediate layer for the current stage. We numerically show that using different reduced order models as inputs of each stage can improve the training and we propose several ways of adding different information into the systems. These methods include mathematical multiscale model reductions and network approaches; but we found that the mathematical approach is a systematical way of decoupling information and gives the best result. We finally verified our training methodology on a time dependent nonlinear problem and a steady state model. 
\end{abstract}

\section{Introduction}
Deep learning has been successfully applied to solve problems with multiscale heterogeneities  \cite{brunton2019methods,cheung2020deep,heinlein2019machine,magiera2020constraint,regazzoni2020machine,vasilyeva2020learning,wang2018multiscale,wang2020reduced,wang2020recurrent,wang2020deep,yeung2020deep,zhang2020learning}. One of the benefits of a deep learning approach is computation time. Once the model is trained, the prediction is a set of multiplication and addition operations; hence the deep learning approach can be used in many areas of multiscale problems. However, there are two major difficulties for researchers to use deep learning in solving multiscale problems.

The first issue is the network design. Researchers now understand the convolutional operation in many computer vision (CV) tasks and there are several effective deep learning models which can be used as the front end in many CV tasks \cite{segnet,dong,resnet,senet,dense,condition,perceptual,long,simonyan2014very,scaling,inception,rethinking}; however, as far as we know, there is no deep learning model nor operation which can be commonly used in physical multiscale problems. Therefore, people are lack of ``state of art" models and need to design a model with the uncertain performance. Even though the network is well designed following some {\it rules of thumb}, we still cannot estimate the potential of the network due to hyper-parameters and an uncertain optimization process.

The second problem is the input and output dimensions. To capture the multiscale nature of the problem, very fine scale information should be used. This usually results in a huge output and input dimension. For example, to solve a time-dependent partial differential equation (PDE), we usually need the physical properties of the PDE such as the heat source or solutions evaluated at the fine scale mesh at all time steps. The input and output dimension will raise the difficulty in the training process. It will also be hard to design the network structure, in particularly the front-end feature extraction (dimension reduction) module which is standard in any network.

To tackle the second problem, instead of learning the fine scale information using high dimensional input, we can learn the reduced order model of the fine problem. We can then easily recover the fine scale model. 
Convergence and accuracy theory have been established for many theoretical multiscale model reduction methods \cite{gmsfem2,snapshot,gmsfem9,chung2010reduced,chung2014generalized,chung2015residual,gmsfem4,chung2018constraint,gmsfem6,chung2016mixed,gmsfem1,gmsfem7,gmsfem8,msfem1}. In this work, we are going to train models aimed at predict the reduced order model of multiscale problems. 

The issue with the output dimension can be solved by predicting the reduced order models; but we still need to deal with the dimension of the input. One idea is to train the deep learning model with the reduced order representation of the fine information; we will develop our method basing on training the multiple reduced order models in multiple stages.


\subsection{Main contributions}

In this work, we propose a multi-stage training strategy for problems with multiscale properties to tackle the issues mentioned above. In the first stage of the training process, a rough prediction is generated. In the following stages, one can iteratively improve the prediction accuracy. Each stage shares an almost identical network structure. 
Since the network structures are fixed for each stage, the proposed method provides an efficient alternative to network design and reduces the work on tuning hyper-parameters of training.

In order to reduce the input dimension in the front-end, a multiscale model reduction is employed. One may choose data-driven model reduction like auto encoder \cite{goodfellow}. 
However, designing and training such networks requires some additional knowledge and works for multiscale problems. In this work, we construct a reduced-order model based on the constraint energy minimizing generalized multiscale finite element method (CEM-GMsFEM) \cite{chung2018constraint}. This provides a simple and reliable approach for reducing dimension in multiscale modeling. If we use the same reduced-order model as input in each stage, we can view this approach as an iterative correction process. 

Furthermore, we find that decoupling the input information and using different reduced order models in different stages are computationally friendly and gives an acceptable (or even better) prediction. 
To understand how information decoupling and multi-stage training work, we consider a linear problem, for example, a linear PDE with a source function. We can solve sub-problems by using decoupled information from the original source; then we combine the sub-solutions. Our proposed deep learning approach is imitating this process. For the general nonlinear problems, the combination can be fulfilled by the nonlinear network through training. 

To summarize, we show by experiments that our training methodology can improve the predictions with the
same reduced order model in each stage. Several enhancements of the proposed methodology
are provided and discussed. One noticeable discovery and upgrade is to decouple the input
information and use different decoupled reduced order models as input to each stage. This
upgrade can save computation cost while providing an even better prediction. Several other
upgrades using a trained deep model are also presented.

To verify the proposed method, we will first verify our proposed methodology on two time-dependent multiscale PDE problems. We design a three-module network structure, which works very well for the time-dependent problems considered in this work. These three modules include dimension reduction, multi-head attention and generator modules. 
With slight modification, we are able to extend our multi-stage framework for steady-state model problem with input being permeability fields. 

\subsection{Outline of the paper}
The rest of the paper is organized as follows. 
In Section \ref{model}, we first present model problems which will be used to test and demonstrate our methodology. 
In Section \ref{deeplearning}, we propose a network architecture of multi-stage learning-based algorithm designed for coupling reduced-order models. 
Numerical experiments related to time-dependent model problems are presented in the Section \ref{experiment}. In the following Section \ref{kappa_model}, we will verify our proposed method on a different multiscale problem. Concluding remarks are drawn in Section \ref{conclusion}. 

\section{Model problems and multiscale model reduction}\label{model}
In this section, we present model problems which will be used to illustrate the proposed deep learning methodology. We will also introduce a multiscale model reduction technique based on constraint energy  minimizing generalized multiscale finite element method (CEM-GMsFEM) developed in \cite{chung2018constraint}. The CEM-GMsFEM will serve as a model reduction method in Section \ref{deeplearning}.

\subsection{Model problem}
Let $\Omega \subset \mathbb{R}^d$ be a bounded domain. 
Let $T >0$ be a fixed time and denote $\mathbb{I} := (0,T]$. We consider the following evolution equation: find $u : \mathbb{I} \times \mathbb{R}^d \to \mathbb{R}^q$ ($q \in \{1, d, ~ d+1 \}$) such that 
\begin{eqnarray}
\mathbf{M} u_t + \mathcal{L}^{\eps} u = f \quad \text{in } \mathbb{I} \times \Omega, \quad u = 0 \quad \text{on } \mathbb{I} \times \partial \Omega, \quad \text{and} \quad 
u = g \quad \text{on } \{ 0 \} \times \Omega.
\label{eqn:model}
\end{eqnarray}
Here, $\mathbf{M} \in \mathbb{R}^{q \times q}$ is a positive definite linear operator, $\mathcal{L}^{\eps}$ denotes a differential operator encoded with multiscale features, and $f$ is a random source function with sufficient regularity. 
Denote by $V(\Omega)$ the appropriate solution space, and $V_0(\Omega) := \{ v \in V(\Omega): v = 0 ~ \text{on} ~ \partial \Omega \}$. The variational formulation of \eqref{eqn:model} reads: find $u(t,\cdot) \in V_0(\Omega)$ such that 
$$ ( \mathbf{M} u_t, v ) _{\Omega} + \left \langle \mathcal{L}^{\eps}u, v \right \rangle_{\Omega} = (f,v)_{\Omega} \quad \tforall  v \in V_0(\Omega),$$
where $\left \langle \cdot , \cdot \right \rangle_{\Omega}$ denotes a specific bilinear form related to the differential operator $\mathcal{L}^{\eps}$ and $(\cdot, \cdot)_{\Omega}$ is the standard $L^2$-inner product. 
Next, we introduce some examples for \eqref{eqn:model}. 
\begin{itemize}
\item {\bf Heat equation.} In this case, it is a parabolic equation with homogeneous Dirichlet boundary condition and we have $q = 1$, $\mathbf{M} = 1$, $\mathcal{L}^\varepsilon (u) = - \nabla \cdot ( \kappa \nabla u)$, $V(\Omega) = H^1(\Omega)$,  $\kappa$ is a heterogeneous space-time (in general) function,
 and $\left \langle \mathcal{L}^\varepsilon (u), v \right \rangle_{\Omega}  = ( \nabla u , \nabla v)_\Omega$.

\item {\bf Richards equation.} In this case, the problem is nonlinear and we have $q = 1$, $\mathbf{M} = 1$, $\mathcal{L}^\varepsilon (u) = - \nabla \cdot (\kappa(u) \nabla u)$, $V(\Omega) = H^1(\Omega)$. Here, $\kappa (u)$ is a heterogeneous space-time function depending on the solution $u$. 
\end{itemize}

\subsection{Model reduction with CEM-GMsFEM}
In this section, we outline the framework of the CEM-GMsFEM, which provides a multiscale model reduction to generate front-end inputs. The procedure of the CEM-GMsFEM can be summarized as follows: (i) construct auxiliary modes; and (ii) construct CEM modes based on auxiliary modes. 

We first introduce the notations of coarse and fine grids. 
Let $\mathcal{T}^H$ be a coarse-grid partition of the domain
$\Omega$ with mesh size $H>0$. Note that the coarse-grid size
is much larger than the size of small-scale heterogeneities in multiscale applications. 
By conducting a conforming
refinement of the coarse mesh, one obtains a fine mesh $\mathcal{T}^h$ with mesh size $h>0$.
Typically, we assume that $0 < h \ll H$, and that the fine-scale mesh $\mathcal{T}^h$
is sufficiently fine to fully resolve the small-scale features;  
while $\mathcal{T}^H$ is a coarse mesh containing many fine-scale features.
Let $N$ and $N_c$ be the number of elements and nodes in coarse grid, respectively. We denote $\{x_i: \, 1\leq i \leq N_c\}$ the set of coarse nodes and $\{K_i: 1 \leq i \leq N\}$ the set of coarse elements. 


\subsubsection*{\underline{Auxiliary modes}}\label{auxiliary}
Let $K_i \in \mathcal{T}^H$ be a coarse element for $i \in \{ 1, \cdots, N\}$. 
We consider the local eigenvalue problem over the coarse element $K_i$ as follows: find $\phi_j^{(i)}\in V(K_i)$ and $\lambda_j^{(i)} \in \mathbb{R}$ such that 
\begin{eqnarray}\label{eq:spectral}
\left \langle \mathcal{L}^\eps \phi_j^{(i)} , v \right \rangle_{K_i} = \lambda_j^{(i)} s_i \left ( \phi_j^{(i)}, v \right ) \quad \text{for all } v \in V(K_i).
\end{eqnarray}
where $s_i(\cdot,\cdot)$ is defined as follows:
$s_i(p,q) := \int_{K_i} \tilde \kappa pq ~ dx$ and $\tilde \kappa := \kappa \sum_{j=1}^{N_c} \abs{\nabla \chi_j}^2.$
Here, $\{ \chi_j \}_{j=1}^{N_c}$ is a set of standard multiscale basis functions satisfying the property of partition of unity. We remark that the definition of $\tilde \kappa$ is motivated by the analysis. 
We arrange the eigenvalues of \eqref{eq:spectral} in ascending order and pick the first $L_i$ eigenfunctions $\{ \phi_j^{(i)} \}_{j=1}^{L_i}$ corresponding to the small eigenvalues in order to to construct the CEM basis functions in next step. 
\subsubsection*{\underline{Construct CEM basis}}
Define $s(\cdot, \cdot) := \sum_{i=1}^N s_i(\cdot,\cdot)$. 
We then construct CEM basis functions; 
the (localized) CEM multiscale basis function $\psi_{j,\text{ms}}^{(i)}$ satisfies the following minimization problem: 
\begin{eqnarray}
\psi_{j,\text{ms}}^{(i)} := \argmin_{\psi \in V_0(K_i^+)}  \left \{  \left \langle \mathcal{L}^\eps \psi, \psi \right \rangle_{K_i^+} \quad \text{s.t. } s(\psi, p_{j'}^{(i')}) = \delta_{jj'} \delta_{ii'} \right \},
\label{eqn:cem_1}
\end{eqnarray}
where $\delta_{ii'}$ is the Kronecker delta function. Here, $K_i^+ := K_{i,\ell}$ is an oversampled region defined as follows:
$$ K_{i,0} := K_i, \quad K_{i,\ell} := \bigcup \{ K \in \mathcal{T}^H: K \cap K_{i,\ell-1} \neq \emptyset \} \quad \text{for } \ell \geq 1$$
and $\ell$ is a user-defined parameter of oversampling. 
Therefore, the multiscale space defined as $V_{\text{ms}} := \text{span} \{ \psi_{j,\text{ms}}^{(i)}: i = 1, \cdots, N, ~ j = 1, \cdots, L_i \}$ provides an accurate (multiscale) model reduction (in the sense of Galerkin projection) to the solution of the problem \eqref{eqn:model}. 
With this set of localized basis functions, the results in \cite{chung2018constraint} show that one can achieve first-order convergence rate (with respect to the coarse-mesh size $H$) independent of the contrast provided that sufficiently large oversampling size is considered. 

\section{Deep learning}\label{deeplearning}
\begin{figure}[H]
\centering
\includegraphics[scale = 0.50]{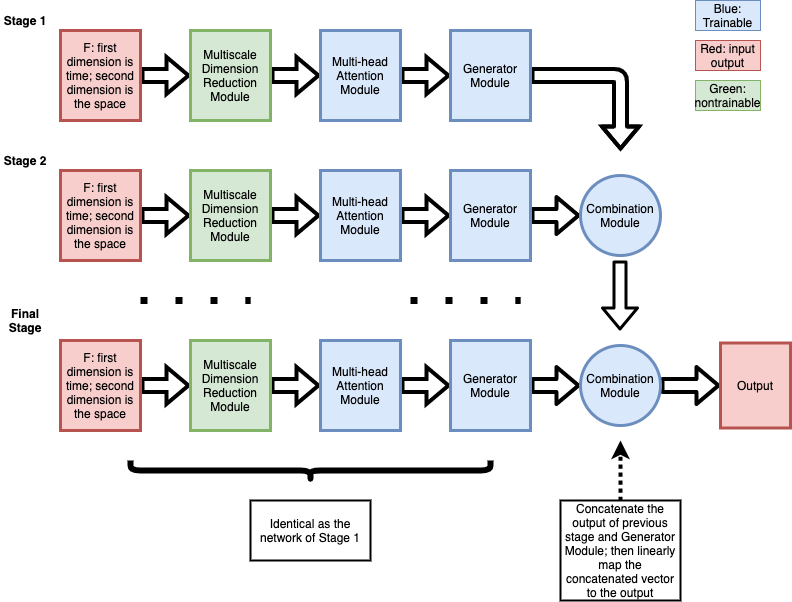}
\caption{Structure and Workflow}
\label{csall}
\end{figure}

\begin{figure}[H]
\centering
\includegraphics[scale = 0.40]{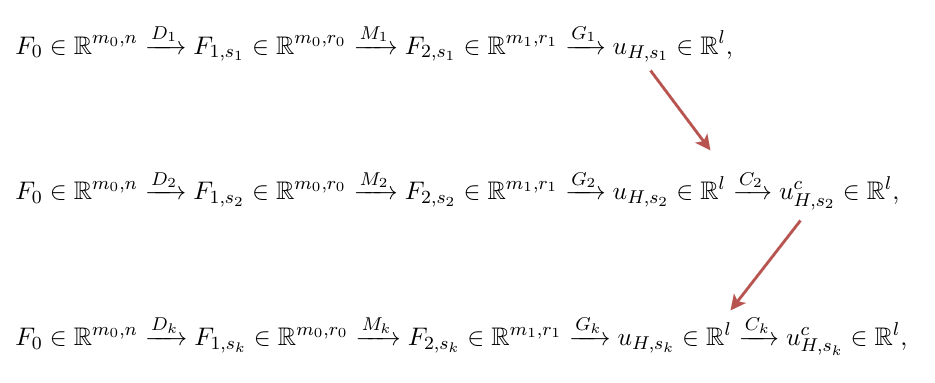}
\caption{Mathematical expression of the structure and workflow; each row represents a stage.}
\label{sall}
\end{figure}

In this section, we will present the multistage training details. 
The idea is to train a rough model in the first stage and make correction based on the prediction in the first stage and additional information. 
The second and the following stages will share the similar structure as that of the first stage but have additional combination module which will combine the prediction of the previous stage and the current stage.
The network structure and workflow are described in Figure \ref{csall}.
The stages share a similar network structure. Except the first stage, the rest of the stages will adapt the intermediate prediction of the previous stage and a reduced order model of the fine information as the input (to trainable part of the network); the input of the first stage is a reduced order model; and the learning target of each stage is the same reduced order model (coarse scale solution). The network structure can be also expressed mathematically as shown in Figure \ref{sall}. 

Note that 
$F_0\in\mathbb{R}^{m_0, n}$ is the input, where $m_0$ is the temporal dimension and $n$ is the spatial dimension. If we are solving the model problem \eqref{eqn:model}, $F_0$ is the source $f(x, t)$ evaluated at all time steps $t_1,...t_{m_0}$ and all fine mesh points $x_1,...,x_n$. $u_{H, s_i}\in\mathbb{R}^l$ is the output of the generator at stage $s_i$ and will be combined with the prediction of the previous stage. 
$u_{H, s_i}^c\in\mathbb{R}^{l}$ is the prediction and target of the learning at stage $i = 1,2,...,k$; it should be noted that $n\gg l$. If we want to solve the model problem \eqref{eqn:model} using CEM-GMsFEM, this can be the coarse scale solution (CEM solution) at the terminal time. Detailed application of the deep model proposed here on the model problem can be seen in Section \ref{model_deep_realization}.

Moreover, we define $D_k: \mathbb{R}^{m_0, n}\rightarrow \mathbb{R}^{m_0, r_0}, M_k: \mathbb{R}^{m_0, r_0}\rightarrow\mathbb{R}^{m_1, r_1}, G_k: \mathbb{R}^{m_1, r_1}\rightarrow\mathbb{R}^{l}, C_k: \mathbb{R}^{l}\times \mathbb{R}^l\rightarrow\mathbb{R}^{l}$ to be Multiscale Dimension Reduction Module, Multi-head Attention Module, Generator Module and Combination Module respectively. 
We will detail each module below.

We consider using the identical network structure in all stages and the motivation is: the network cannot be well trained due to the limitations of optimizer and tons of hyper-parameters.
Multi-stage training then provides a different initialization and some additional information (output of the previous stage) to the network and can hence improve the training of the network. 

We also propose several ways to upgrade the fundamental approach; all the upgrade options adopt the identical network structures as the fundamental method; but we will decouple the input information.
The motivation is: how to correct the prediction using the same network? Our idea is to use different information and combine the predictions. 

To understand how information decoupling and multi-stage training work, we consider a linear problem, for example, a linear PDE with a source function. We can solve it using decoupled information from the original source and linearly combine the sub-solutions. Our methodology mimics this process in the Combination Module through training. We remark that this methodology is also available to nonlinear problems.


\subsection{Multiscale Dimension Reduction Module}\label{multiscale_dim_red}
In this section, we introduce the multiscale dimension reduction module of the network.
In order to capture multiscale features in multiscale problems, one has to involve an extremely high degrees of freedom and it leads to a large input in the network. 
The multiscale dimension reduction module is targeted to reduce the dimension of the input. We remark that this module is independent of the training process. Hence, it will benefit the front-end (features extraction) design and the training of the network. Mathematically speaking, the multiscale dimension reduction module $D_k$ has is a linear transformation defined below: 
$$
D_k: \mathbb{R}^{m_0, n}\rightarrow \mathbb{R}^{m_0, r_0}.
$$
where $r_0 \in \mathbb{N}$ and $r_0<n$. In particular, the spatial dimension will be reduced through the module. Modern networks usually can be decomposed into a dimension reduction end and a functional end. Since the input data are of multiscale nature (either in the sense of data or mathematics), we design an untrained module which will reduce the dimension of the raw input. Intuitively, the reduced dimension data of each stage should be different since we use the identical network structure. We hence need a systematical way to prepare the reduced dimension data.
We propose three candidates for the choices of reduced order models.

\subsubsection{Mathematical Approach}\label{idea1}
The coarse information can be obtained by mathematical multiscale model reduction (MMMR). One can use CEM-GMsFEM introduced above to form the multiscale dimension reduction module. We denote $R = [\psi_{j, ms}^{i}]\in\mathbb{R}^{n, r_0}$, where $j = 1,..., L_i $ and $i = 1,...,N$. Then, we define $D_k$ such that 
$$D_k(F_0) = F_0R\in \mathbb{R}^{r_0}$$
It should be noted that temporal dimension can also be reduced if the problem has multiscale nature in time. One of the benefits of the MMMR is that it is a systematical way and we can control the input of each stage so that the inputs are orthogonal to each other. In fact, our experiments show that the MMMR is the key in the success of the multi-stage training.

\subsubsection{Raw Data Approach}\label{idea2}
Instead of using certain mathematical model to get the coarse information, we can acquire it by performing max-pooling on the fine grid information. Max-pooling operator is a local operation which will take the max value in a local neighborhood. 
\textcolor{black}{Due to the multiscale nature of our problem, max-pool is a great candidate for extracting information from the fine grid.}
The disadvantage is also apparent: although we can control the hyper-parameters (pooling size, stride etc.) in the pooling operations, the pooled data with different hyper-parameters are similar. We hence consider using the pooling as a complement information to our network.

\subsubsection{Data-Driven Approach}
We can train an unsupervised deep learning model to achieve model reduction; and later use the intermediate layers as the reduced order model. \textcolor{black}{One successful deep model reduction structure should be auto-encoder and under some assumptions, people can show that linear network is equivalent to the principle component analysis.}
Let us define
$$\mathcal{A}: u_h\rightarrow u_h\in\mathbb{R}^n$$
as the feature extraction network stacked of multiple layers as standard (see Figure (\ref{ae}) for the illustration of the network), where $u_h$ is the fine scale solution and is the target of the prediction; if we denote $\mathcal{A}_i$ as the  output of $ith$ layer of network $\mathcal{A}$, the dimension reduction can be expressed as:
$$
D_k(F_0) = \mathcal{A}_i.
$$
Compared to the mathematical multiscale model reduction, trained model may give better prediction (see the experiments in Section \ref{exp4}); however,  there are some issues with the trained network.
\begin{enumerate}
\item Designing the model reduction network. We are lack of rules of thumbs in the multiscale problems. One of the purposes of this work is to reduce the workload in designing the network. In fact, we designed several networks in our research but the performances are very diverse. We created a network which gives the best result among all experiments; however, most networks we designed fail to give us a good model reduction comparing to the other two approaches.
\item Picking up the intermediate layers. Even we have a decent model reduction network, it is not easy to choose layers which can be served as the reduced order models. Although it is clear in the computer vision area \cite{condition, perceptual} about layers information, there is no common understanding in the multiscale applications of the deep neural networks.
\end{enumerate}

The above two mentioned points are interesting and we will study them in the future research.

\subsection{Multi-head Attention Module}\label{mha}
In this section, we briefly introduce the multi-head attention module in the proposed network structure. In particular, the technique of multi-head attention \cite{vaswani2017attention} will be used to further reduce the spatial and temporal dimensions of the problem. 
Multi-head attention is originally designed to deal with the natural language processing problems. It can be shown that this technique encodes the long time dependency of the sequence and reduce the dimension of the sequence. For example, many words have multiple meanings. In order to accurately encode the meaning of a word (into a vector), we need to look at the whole sentence; that is, we need to check the connection between this word and the other words in a sentence. In short, multi-head attention is designed to return a more concise and accurate representation of this sentence. 

In the proposed time dependent problem, the vector $F_0\in\mathbb{R}^{m_, n}$  is the input source $f(x, t)$ evaluated at all time steps and fine mesh points; so it is is natural to consider the connections and the joint contributions of the source at different time steps. Here, we use the technique of multi-head attention to give a shorter and accurate representation of the whole process by including rich information of the long time dependency of the sources at different time. For example, the effect of the $f(x, t_0)$ on $f(x, t_{m_0})$ will be encoded into the output vector. 
We define the multi-head attention module $M_k$ such that
$$
M_k: \mathbb{R}^{m_0, r_0}\rightarrow\mathbb{R}^{m_1, r_1}
$$
where $m_0>m_1$  and $r_0>r_1$. The pair $(m_0, m_1)$ represents the original space-time dimension before reduction while $(m_1, r_1)$ refers to the reduced dimension transformed by the multi-head attention.  
The output of this module will be denoted as $F_2^{m_1, r_1}$, where $m_1$ is the reduced temporal dimension and $r_1$ is the spatial dimension.
We use the classical transformer with $6$ heads and $1$ layer in the simulations below. 
This is the second dimension reduction and the output will be used as the input to the back end generator module. It should be noted that if the problem has no time dependency, this multi-head attention module will be removed.

\subsection{Generator Module}
The first two modules are target to conduct the dimension reduction and feature extraction. The Generator Module, however, is designed to generate the target solution $u_{H, s_k}$ from the low dimensional representation of the input; it will enlarge the dimension of the feature extracted from the front end and map it to the final prediction.

In our network, we use one layer fully connected network to perform the generation and our experiments show that if the front-end is powerful enough to encode the information of the input time dependent process, the structure of the generator can be fairly simple. The output of the generator $G_k$ of the first stage will be the prediction; however, for the later stage, the output of the generator will be combined with the prediction of the previous stage in the combination module which will be explained later.

\subsection{Combination Module}
Combination module is targeted to combine the prediction of previous stage $u_{H, s_{k-1}}^c$ and the output of the current Generator $u_{H, s_k}$. We concatenate $u_{H, s_{k-1}}^c$ and $u_{H, s_k}$ and map it to the final prediction; i.e.,
$$
C_k: \mathbb{R}^{l}\times \mathbb{R}^l\rightarrow\mathbb{R}^{l}.
$$
The motivation is that the prediction of the current stage will correct the previous prediction in an iterative way (due to the same network structure). This can also be taken as a trainable ensemble process; however, it should be noted that we use the same model but different input information in each model, which is different from the classical ensemble. 

The design of the combination module and the proposed multi-stage training using decoupled information can also be motivated mathematically.
If we solve a linear problem and use the different inputs in each sub-problem, the Combination Module can be taken as a linear combination of the sub-solutions. For the general nonlinear problems, the Combination Module is an extension of the linear combination and can be fulfilled through the training process.



\section{Numerical Experiments}{\label{experiment}}
In this section, we test our multi-stage learning algorithm with several linear and nonlinear problems that are representative in multiscale modeling. We present some experiment results of the proposed learning strategy. 
We first demonstrate an initial version in Section \ref{exp1}; several upgrades which are based on the idea of decoupling information will be provided in Sections \ref{exp1_1}, \ref{exp3}, and \ref{exp4}. A comment on the efficiency of the upgraded version is given in Section \ref{comment_eff}. We also evaluate the error of the multiscale method (using learnt coarse scale model) and show the results in Section \ref{err_ms}.

\subsection{Experiment Settings}
We set $\Omega = (0,1)^2$ and $T = \pi$. Time step is set such that the total number of time steps is equal to $m_0 = 31$. We partition the spatial domain into $10 \times 10$ equal square elements to form the coarse grid $\mathcal{T}^H$ and thus $H = 0.1$. Similarly, we set $h = 0.01$ to form the fine grid. We set $L_i = 3$ to form the multiscale space $V_{\text{ms}}$. 
In this case, we have $\mathcal{N} := \text{dim}(V_{\text{ms}}) = 300$ and we make use of this multiscale space to achieve multiscale model reduction. We use $1200$ samples to train the network and test on other $400$ samples. 

\begin{figure}[h]
\centering
\mbox{
\includegraphics[scale = 0.5]{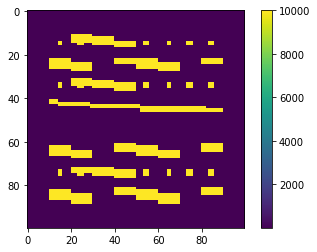}
\quad 
\includegraphics[scale = 0.5]{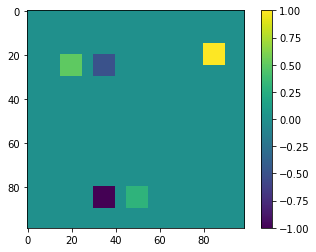}
}
\caption{Permeability field $\kappa$ (left); the support of the source function $f$ (right). }
\label{kappa}
\end{figure}

We first consider a linear problem when $q = 1$, $\mathbf{M} = 1$, and $\mathcal{L}^\eps u:= - \nabla \cdot (\kappa (x) \nabla u)$ in \eqref{eqn:model}. 
The permeability field $\kappa : \Omega \to \mathbb{R}$ is defined as shown in Figure \ref{kappa} (left).

The linear equation to be considered reads: find $u \in \mathbb{I} \times \Omega \to \mathbb{R}$ such that 
\begin{align}
u_t = \nabla\cdot (\kappa \nabla  u) +f \quad \text{in } \mathbb{I} \times \Omega, \quad u = 0 \quad \text{on } \mathbb{I} \times \partial \Omega, \quad \text{and} \quad 
u = 0 \quad \text{on } \{ 0 \} \times \Omega.
\label{eqn:linear-model}
\end{align}
The source function $f$ is equipped with randomness and its support is depicted in Figure \ref{kappa} (right). 
In particular, the configuration of $f$ is listed as follows: 
\begin{eqnarray}
f(t,x) = \left \{  
\begin{array}{ll}
e^{1+ \cos (t)} + \xi_1 & \text{for } x \in [0.15, 0.25] \times [0.8, 0.9], \\
\\
- e^{1+ \cos (t)} + \xi_2 & \text{for } x \in [0.8, 0.9] \times [0.3, 0.4], \\
\\
6 \cos (t) + \xi_3 & \text{for } x \in [0.2, 0.3] \times [0.15, 0.25], \\
\\
\displaystyle{\frac{1}{2} t^2 + 0.5} + \xi_4 & \text{for } x \in [0.8, 0.9] \times [0.45, 0.55], \\
\\
\displaystyle{-\frac{1}{2} (\pi - t)^2} + \xi_5 & \text{for } x \in [0.2, 0.3] \times [0.3, 0.4], \\
\\
0  & \text{otherwise},
\end{array}
\right .
\label{eqn:source}
\end{eqnarray}
where $\{ \xi_i \}_{i=1}^5$ is a set of independent and identically distributed random variables with uniform distribution $U(-0.5, 0.5)$. We refer \eqref{eqn:linear-model} as the {\it linear case} in the simulation below. 

We also consider a nonlinear problem with $q = 1$, $\mathbf{M} = 1$, and $\mathcal{L}^\varepsilon (u) := - \nabla \cdot (a(x,u) \nabla u)$ in \eqref{eqn:model}, where $a(x,u) := \kappa(x) \exp (\gamma u)$ is a nonlinear permeability field with $\kappa$ depicted in Figure \ref{kappa} (left) and $\gamma$ being a positive constant. The nonlinear model problem to be solved reads: 
\begin{align}
u_t = \nabla\cdot (a(x,u) \nabla  u) +f \quad \text{in } \mathbb{I} \times \Omega, \quad u = 0 \quad \text{on } \mathbb{I} \times \partial \Omega, \quad \text{and} \quad 
u = 0 \quad \text{on } \{ 0 \} \times \Omega.
\label{eqn:nonlinear-model}
\end{align}
The source function $f$ is defined in \eqref{eqn:source}. We set $\gamma = 20$ in the experiments and refer  \eqref{eqn:nonlinear-model} as the {\it nonlinear case} in the simulation. 

\subsection{Learning Details}\label{model_deep_realization}
We solve the linear case \eqref{eqn:linear-model} or the nonlinear case \eqref{eqn:nonlinear-model} combining with the proposed network structure and learning strategy. 
Instead of predicting the fine grid solution, we predict a coarse mesh solution $u_H$ to the problem. 


Suppose that $\{ x_j \}_{j=1}^n$ is the collection of all fine-grid nodes in $\mathcal{T}^h$; and the temporal domain $\mathbb{I} = (0,T]$ is divided into $m_0 \in \mathbb{N}^+$ pieces. We denote $F_0\in\mathbb{R}^{m_0, n}$ as $(F_0)_{ij} := f(t_i, x_j)$ for $x_j$ being fine-grid nodes and $t_i$ being discrete time level; the input $F_0$ is the projection of the source at different time onto the fine scale mesh. Then, the target of all experiments can be summarized as

$$
F_0\xrightarrow{    \text{Deep Neural Network with the Proposed Methodology}    } u_H(t_{m_0}) := (R^tR)^{-1}R^t u_h(t_{m_0})\in\mathbb{R}^N. 
$$
Here, $R\in\mathbb{R}^{n, \mathcal{N}}$ is the matrix collecting all basis functions in $V_{\text{ms}}$ with representation on the fine grid;
the learning target is the coarse-scale model and is the projection of the fine scale solution onto the CEM space $V_{ms}$. 
The function $u_h(t_{m_0})$ denotes the fine-grid solution evaluated at the terminal time $t_{m_0}$. 
In the first stage, we solve the following minimization problem: 

$$
\min_{G_1, M_1} \norm{u_H(t_{m_0}) - G_1(M_1(D_1(F_0)))}_1,
$$
where $\norm{v}_1 := \sum_{i = 1}^{N} |v_i| $ for $v\in\mathbb{R}^N$.
In the following stage $k$ with $k \geq 2$, we solve 
$$
\min_{G_k, ~M_k, ~C_k} \norm{u_H(t_{m_0})- C_k(u_{H, s_{k-1}}^c, u_{H, s_{k-1}} )}_1
$$
with $u_{H, s_{k-1}} = G_k(M_k(D_k(F_0))) \in\mathbb{R}^N$ and $u_{H, s_{k-1}}^c \in\mathbb{R}^N$ being  the prediction generated from the previous stage. 

The mappings $D_k$ are Dimension Reduction Modules. We will test and compare all dimension reduction methods described in Section \ref{multiscale_dim_red} in details in the experiments below. For example, our baseline model is: $D_k(F_0) = R^t F_0$ for all $k$, i.e., the projection of input onto the entire basis in all stages; this is demonstrated in the first experiment. 

We remark that the performance of the method will be measured by relative vector $l_2$ norm which is defined as follow:
\begin{eqnarray}
    \text{relative $l_2$ error} = \frac{ \|u_{H}-u_{H, learnt}\|_2 }{\|u_H\|_2},\label{relatve_l2}
\end{eqnarray}
where $u_H := (R^tR)^{-1}R^t u_h$ is the exact coarse-grid solution and $u_{H, learnt}$ is the coarse-grid prediction generated by the proposed learning algorithm. Here, $\norm{\cdot}_2$ is the standard Euclidean norm for vectors. One possible way to test the performance of our method is to compute the mean of the training data and then make predictions using the mean. We calculate the relative error using the mean for both linear and nonlinear cases. The relative $l_2$ errors are $0.13567$ for the linear problem and $0.17638$ for the nonlinear problem.

\subsection{Experiment 1: multi-stage correction}\label{exp1}

In this set of experiments, we are going to show the proposed multi-stage training can improve the prediction even if we use same input information in each stage. This is baseline of the method and can be improved by using decoupled information which will be showed in later sections. The network structure of each stage is fixed and the input will be projection of the source onto the CEM multiscale basis. To be more specific, we use CEM-GMsFEM with $H = 0.1$ and $L_i = 3$ for all $i$ and all stages; that is, the Module $D_k$ for all stage $k$, are the projection on all CEM basis functions. 

We will present the mean relative error as defined in \eqref{relatve_l2} for each stage. We modify the structure of the network and compute the error for each one. Networks differ from each other by $(m_1, r_1)$ which are defined in Section \ref{mha} and are important hyper-parameters regarding multi-head attention. This network variation convention will be used in all experiments in the later sections.
Tables \ref{table1_1} and \ref{table1_1nl} show the results for linear and nonlinear cases respectively. Columns $2$ to $4$ represent the error of each stage.
\begin{table}[H]
\centering
\begin{tabular}{||c c c c||} 
\hline
$(m_1, r_1)$ & Stage 1 Error & Stage 2 Error & Stage 3 Error \\ [0.5ex] 
\hline
(30, 10)  & 0.13995 & 0.09439 & 0.08148 \\ [0.5ex]
\hline
(25, 12)  & 0.12647 & 0.07955 & 0.07399 \\ [0.5ex]
\hline
(20, 15)  & 0.11721 & 0.07453 & 0.06826 \\ [0.5ex]
\hline
(15, 20)  & 0.11405 & 0.08231 & 0.06647 \\ [0.5ex]
\hline
(10, 30)  & 0.12135 & 0.07069 & 0.06670 \\[0.5ex]
\hline
\end{tabular}
\caption{(Linear) Entire $3$ CEM basis ($H = 0.1$ and  $L_i =3$) are used in all stages.}
\label{table1_1}
\end{table}

\begin{table}[H]
\centering
\begin{tabular}{||c c c c||} 
\hline
$(m_1, r_1)$ & Stage 1 Error & Stage 2 Error & Stage 3 Error \\ [0.5ex] 
\hline
(30, 10)  & 0.10287 & 0.10118 & 0.09887 \\ [0.5ex]
\hline
(25, 12)  & 0.11385 & 0.09747 & 0.09201 \\ [0.5ex]
\hline
(20, 15)  & 0.18628 & 0.09920 & 0.09267 \\ [0.5ex]
\hline
(15, 20)  & 0.12656 & 0.09957 & 0.09043 \\ [0.5ex]
\hline
(10, 30)  & 0.10555 & 0.09413 & 0.09262 \\[0.5ex]
\hline
\end{tabular}
\caption{(Nonlinear) Entire $3$ CEM basis ($H = 0.1$ and $L_i =3$) are used in all stages.}
\label{table1_1nl}
\end{table}
It should be noted that we train the first stage with so many epochs that the loss decays very slowly at the end of the training. 
We can compare the stage 1 and stage 3 error and in all cases (different network structures), the predictions can be improved by the multistage training. 

Intuitively, the network is re-initialized with additional input and this may help the training. We can observe that the improvement from stage 2 to stage 3 is marginal. This probably because the iterative correction process almost converges to its limit and hence has small improvement.

\subsubsection{A Stable Upgrade}\label{exp1_1}
In this set of experiments, we will present an upgrade of the initial multistage training, i.e., we are going to decouple the input information and use different reduced order models in each stage.
This is more natural as we use different information but same network to correct the predictions.
Compared to our initial idea, the input size is reduced and hence we have smaller network; this will save GPU memory and own less trainable parameters. However, the accuracy is not compromising and we observed even better prediction when we use the same setting of the training as the first set of experiments.

We will choose different multiscale basis progressively in three stages and we use CEM-GMsFEM with $H = 0.1$ and $L_i = 3$ for all $i$.
The dimension reduction module of the first stage is trained with source $f(x, t)$ projected onto the first basis in each coarse neighborhood; the second and third stages are then trained with the projection onto the second and third local basis progressively. Network structures for each stage are identical except the combination module in the second and third stage.

Similarly as before, we present the mean relative error for each stage 
and test on different networks following the same convention as before.
See Table \ref{table1} for the linear problem and \ref{table1nl} for the nonlinear problem. 
\begin{table}[H]
\centering
\begin{tabular}{||c c c c||} 
\hline
$(m_1, r_1)$ & Stage 1 Error & Stage 2 Error & Stage 3 Error \\ [0.5ex] 
\hline
(30, 10)  & 0.13572 & 0.05322 & 0.04280 \\ [0.5ex]
\hline
(25, 12)  & 0.13282 & 0.07162 & 0.06090 \\ [0.5ex]
\hline
(20, 15)  & 0.13932 & 0.06584 & 0.05996 \\ [0.5ex]
\hline
(15, 20)  & 0.12529 & 0.06231 & 0.04656 \\ [0.5ex]
\hline
(10, 30)  & 0.13904 & 0.05116 & 0.04182 \\[0.5ex]  
\hline
\end{tabular}
\caption{(Linear) $1^{st}$, $2^{nd}$ and $3^{rd}$ CEM basis ($H = 0.1$ and $L_i =3$) are used one in each stage. }
\label{table1}
\end{table}

\begin{table}[H]
\centering
\begin{tabular}{||c c c c||} 
\hline
$(m_1, r_1)$ & Stage 1 Error & Stage 2 Error & Stage 3 Error \\ [0.5ex] 
\hline
(30, 10)  & 0.13322 & 0.08701 & 0.08439 \\ [0.5ex]
\hline
(25, 12)  & 0.11515 & 0.08538 & 0.08162 \\ [0.5ex]
\hline
(20, 15)  & 0.09814 & 0.08788 & 0.08210 \\ [0.5ex]
\hline
(15, 20)  & 0.12591 & 0.08538 & 0.08054 \\ [0.5ex]
\hline
(10, 30)  & 0.29105 & 0.09888 & 0.07815 \\[0.5ex]  
\hline
\end{tabular}
\caption{ (Nonlinear) $1^{st}$, $2^{nd}$ and $3^{rd}$ CEM basis ($H = 0.1$ and $L_i =3$) are used one in each stage. }
\label{table1nl}
\end{table}

\subsubsection{A Comment On The Improved Efficiency}\label{comment_eff}
In the last set of the experiments, we have demonstrated the accuracy of the upgrade of the multi-stage training. In this section, we briefly show the improved efficiency of the upgrade version. The difference is the input of each stage is decoupled and the upgraded version has smaller input and hence smaller network and memory. In Table \ref{nb_para}, the number of trainable parameters is presented for the initial version and upgraded version; this indicates the complexity of the network.
\begin{table}[H]
\centering
\begin{tabular}{||c c c||} 
\hline
$(m_1, r_1)$ & Without Decouple & Decoupled \\ [0.5ex] 
\hline
(30, 10)  & 235390 & 199390  \\ [0.5ex]
\hline
(25, 12)  & 246492 & 203292  \\ [0.5ex]
\hline
(20, 15)  & 263235 & 209235  \\ [0.5ex]
\hline
(15, 20)  & 291380 & 219380  \\ [0.5ex]
\hline
(10, 30)  & 348570 & 240570  \\[0.5ex]  
\hline
\end{tabular}
\caption{ Number of trainable parameters associated with the network. Column 2 is the original version which are trained using 3 basis; Column 3 is the upgraded version with decoupled input.}
\label{nb_para}
\end{table}

\subsection{Experiment 2: Further Upgraded Option}\label{exp3}
As we have discussed in Section \ref{idea2}, we can use the pooled raw data as input.
This inspires us using the pooled raw data as a complement to the system; i.e., using the pooled raw data directly in the last stage of the training.
To be more specific, $D_k$ ($k = 1, 2$) are the input source projected onto the first and second multiscale basis; 
while in the last stage, $D_3$ will be the pooling of the input data. 

The pooling size is $10$ and stride is $10$ in our experiments. This setting will generate the inputs of the same size as the 1 CEM basis projection;
however, the pooling has no trainable parameter and compared to the mathematical approach, this upgrade can save time computing one CEM base and computing the projection.
Therefore this upgrade will be more efficient when compared to the last experiment \ref{exp1_1}. The results are shown in Tables \ref{table3} and \ref{table3nl}. 

\begin{table}[H]
\centering
\begin{tabular}{||c c c c||} 
\hline
$(m_1, r_1)$ & Stage 1 Error & Stage 2 Error & Stage 3 Error \\ [0.5ex] 
\hline
(30, 10)  & 0.13572 & 0.05322 & 0.04039  \\ [0.5ex]
\hline
(25, 12)  & 0.13282 & 0.07162 & 0.05452  \\ [0.5ex]
\hline
(20, 15)  & 0.13932 & 0.06584 & 0.05828  \\ [0.5ex]
\hline
(15, 20)  & 0.12529 & 0.06231 & 0.03948 \\ [0.5ex]
\hline
(10, 30)  & 0.13904 & 0.05116 & 0.03433 \\[0.5ex]
\hline
\end{tabular}
\caption{(Linear) $1^{st}$ and $2^{nd}$ stages are trained with the projection onto the first and second multiscale basis; the last stage is trained with pooled raw data with pool size 10 and stride 10.}
\label{table3}
\end{table}

\begin{table}[H]
\centering
\begin{tabular}{||c c c c||} 
\hline
$(m_1, r_1)$ & Stage 1 Error & Stage 2 Error & Stage 3 Error \\ [0.5ex] 
\hline
(30, 10)  & 0.13322 & 0.08701 & 0.08597  \\ [0.5ex]
\hline
(25, 12)  & 0.11515 & 0.08538 & 0.08318  \\ [0.5ex]
\hline
(20, 15)  & 0.09814 & 0.08788 & 0.08563  \\ [0.5ex]
\hline
(15, 20)  & 0.12591 & 0.08538 & 0.07697 \\ [0.5ex]
\hline
(10, 30)  & 0.29105 & 0.09888 & 0.08736 \\[0.5ex]
\hline
\end{tabular}
\caption{(Nonlinear) $1^{st}$ and $2^{nd}$ stages are trained with the projection onto the first and second multiscale basis; the last stage is trained with pooled raw data with pool size 10 and stride 10.}
\label{table3nl}
\end{table}

Our experiments show that raw data works well for the linear case, an further improvement can be observed; however, for the nonlinear case, pooled raw data does not give us a stable improvement when compared to the last experiment in Section \ref{exp1_1}.

\subsubsection{Error of the predicted coarse solutions}\label{err_ms}
We also compute the relative $L_2$ error of the predicted coarse solution and compare these errors with the target coarse solution. The relative $L_2$ error is defined as follows: 
\begin{eqnarray}
\text{Relative $L_2$ error} = \frac{\|R u_{H, learnt} - u_h\|}{\|u_h\|},
\end{eqnarray}
where $u_h$ is the fine grid solution and $\|.\|$ is the $L_2$ norm associated to the model problem \eqref{eqn:model}. The errors can be seen in Tables \ref{li_table3_1} and \ref{table3_1}.

\begin{table}[H]
\centering
\begin{tabular}{||c c c ||} 
\hline
$(m_1, r_1)$ &  Learnt &  Computed  \\ [0.5ex] 
\hline
(30, 10)  & 0.00357 & 0.00018  \\ [0.5ex]
\hline
(25, 12)  & 0.00257 & 0.00018  \\ [0.5ex]
\hline
(20, 15)  & 0.00493 & 0.00018  \\ [0.5ex]
\hline
(15, 20)  & 0.00311 & 0.00018  \\ [0.5ex]
\hline
(10, 30)  & 0.00336 & 0.00018 \\[0.5ex]  
\hline
\end{tabular}
\caption{(Linear) Relative $L_2$ error of the learnt and target coarse scale solution. Columns with "Learnt": coarse solutions predicted by the network. Columns with "Computed": coarse solutions computed theoretically (learning target).}
\label{li_table3_1}
\end{table}

\begin{table}[H]
\centering
\begin{tabular}{||c c c||} 
\hline
$(m_1, r_1)$ &  Learnt &  Computed  \\ [0.5ex] 
\hline
(30, 10) & 0.00807 & 0.00015  \\ [0.5ex]
\hline
(25, 12)  & 0.00767 & 0.00015  \\ [0.5ex]
\hline
(20, 15)  & 0.00756 & 0.00015   \\ [0.5ex]
\hline
(15, 20)  & 0.00626 & 0.00015  \\ [0.5ex]
\hline
(10, 30)  & 0.00830 & 0.00015  \\[0.5ex]  
\hline
\end{tabular}
\caption{ (Nonlinear) Relative $L_2$ error of the learnt and target coarse scale solution. Columns with "Learnt": coarse solutions predicted by the network. Columns with "Computed": coarse solutions computed theoretically (learning target).}
\label{table3_1}
\end{table}

We can observe that the relative $L_2$ error of the predicted coarse scale model is small. One can learn a reduced-order model instead of learning a fine model of much higher dimension.

\subsection{Experiment 3: an unstable upgrade}\label{exp4}
\begin{figure}[h]
\centering
\includegraphics[scale = 0.5]{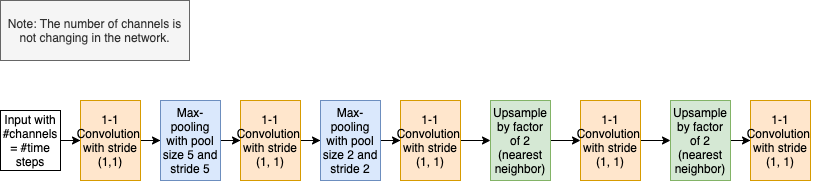}
\caption{Structure of the trained reduced order model network. The first $1-1$ convolution layer is used as the reduced order model.}
\label{ae}
\end{figure}

For all experiments above, we use either mathematical model reduction or pooling; we have seen the importance of the mathematical skills and the best result (in terms of efficiency and for linear case, also the accuracy) is obtained by the combination of 2 information approaches. No training is included in acquiring the low dimension representation; but there are some unsupervised deep learning techniques which can provide a good data representation. This motivates us training a dimension reduction network and the low dimension representation of this network can then be served as the reduced order model.  The network structure can be seen in Figure \ref{ae}. We essentially use the 1 by 1 convolution and max pooling. 

We conduct a two-stage training; $D_1(F_0)$ is an intermediate layer of a pre-trained network. In the network we demonstrate in Figure \ref{ae}, we use the output of the second 1 by 1 convolution layer since this layer has the same dimension of the first stage as the previous experiments. The second stage is trained with the pooled raw data.
As shown in Tables \ref{table6} and \ref{table6nl}, this two-stage model has the best result among all experiments in the previous sections; but we do want to comment that, designing the feature extraction network is not easy; even we can follow some rule of thumbs to design a reasonable network, we have problems picking up a good layer as the reduced order model. In fact, most networks we designed cannot outperform the previous models using mathematical model reduction. There is no guarantee that the feature extraction network you design can really perform a model reduction which can be used in the prediction; we hence call this {\it unstable upgrade}. 

\begin{table}[H]
\centering
\begin{tabular}{||c c c||} 
\hline
$(m_1, r_1)$ & Stage 1 Error & Stage 2 Error \\ [0.5ex] 
\hline
(30, 10)  & 0.11314 & 0.04353  \\ [0.5ex]
\hline
(25, 12)  & 0.10548 & 0.05028  \\ [0.5ex]
\hline
(20, 15)  & 0.11547 & 0.03588  \\ [0.5ex]
\hline
(15, 20)  & 0.10493 & 0.03519  \\ [0.5ex]
\hline
(10, 30)  & 0.11106 & 0.03246 \\[0.5ex]
\hline
\end{tabular}
\caption{(Linear) $1^{st}$ stage is trained with an intermediate layer of a model reduction net; $2^{nd}$ stage is trained with pooled raw data with pool size 10 and stride 10.}
\label{table6}
\end{table}

\begin{table}[H]
\centering
\begin{tabular}{||c c c||} 
\hline
$(m_1, r_1)$ & Stage 1 Error & Stage 2 Error \\ [0.5ex] 
\hline
(30, 10)  & 0.13078 & 0.08360  \\ [0.5ex]
\hline
(25, 12)  & 0.11427 & 0.07887  \\ [0.5ex]
\hline
(20, 15)  & 0.11601 & 0.09122  \\ [0.5ex]
\hline
(15, 20)  & 0.11544 & 0.07940  \\ [0.5ex]
\hline
(10, 30)  & 0.12812 & 0.07343 \\[0.5ex]
\hline
\end{tabular}
\caption{ (Nonlinear) $1^{st}$ stage is trained with an intermediate layer of a model reduction net; $2^{nd}$ stage is trained with pooled raw data with pool size 10 and stride 10.}
\label{table6nl}
\end{table}

We hence conclude that multiple stages training methodology we proposed is a great way to improve the training for the fixed network. We first need to decouple the input information and then train stage by stage. This can save your time in tuning hyper-parameters and designing the network structure. Also, we find that the mathematical model reduction is the key of the success and is the best approach to adding information in the framework of the proposed training methodology.

\section{Applications On Steady-State Model}
\label{kappa_model}
In this section, we present an application of the proposed multi-stage learning algorithm for the diffusion problem with input being permeability field. The previous examples in Section \ref{experiment} are designed with initial input being the external source terms coupled with specific multiscale model reduction. We show here that our proposed strategy can be employed to problems with more complicated input format. 

\subsection{Problem setting}
Consider the following steady-state diffusion problem in a computational domain $\Omega \subset \mathbb{R}^2$: find $u : \Omega \to \mathbb{R}$ such that 
\begin{eqnarray}
- \nabla \cdot (\kappa (x) \nabla u)  = f \quad \text{in } \Omega, \quad u = 0 \quad \text{on } \partial \Omega. 
\label{eqn:steady}
\end{eqnarray}
In this case, the permeability field $\kappa(x) = \kappa(x_1, x_2)$ is assumed to be composed with several different components in the following form $\kappa = \sum_{j} \kappa^j$. 

We employ the proposed learning strategy on the problem \eqref{eqn:steady}. We assume that $\mathcal{T}^h$ and $\mathcal{T}^H$ are the fine and coarse partitions of the domain, with $\mathcal{T}^h$ being a refinement of the coarse grid. 
Under this setting, the permeability field $\kappa$ depends on the random parameters $p_i's$. 
As a result, the multiscale basis functions also depend on the parameters. We obtain multiscale features of the permeability field by taking the average of each scale in each local coarse element. In particular, given a sample of parameter $p_j$, we define 
$$
f_i^j := \frac{1}{N_e^i} \sum_{k=1}^{N_e^i} \kappa^j(x_k, p_j)\chi_{i} (x_k), \quad \text{for } i = 1, \cdots, N,
$$
where $x_k$ is a fine-grid node in each coarse element $K_i$, $\chi_i$ is the characteristic function of the coarse element $K_i$, $N_e^i$ is the number of fine elements in $K_i$, and $N$ is the total number of coarse element.  
The vector $( f_i^j )_{i=1}^N \in\mathbb{R}^N$ will then be used as an input in different stages when we apply the multi-stage strategy to this problem. 
The training target $u_{H,i}$ will be of the form 
$$
u_{H, i} := \frac{1}{N_e^i} \sum_{k=1}^{N_e^i} u_f(x_k) \chi_{i} (x_k), 
$$
where $u_f$ is a reference solution defined on the fine scale. That is the average of the fine solution in each local coarse element. 

In the following experiments, we set $\Omega = (0,1)^2$ and $H = 0.1$. We divide the domain into $10 \times 10$ equal square elements. Thus, we have $N = 100$. The learning process is summarized as follows: 
$$
(f_i^j)_{i=1}^N \xrightarrow{    \text{Deep Neural Network with the Proposed Methodology}} u_H := (u_{H, i})_{i=1}^N,
$$
where $j$ indicates the multiscale feature $\kappa^j(x, p_j)$. We set $\kappa = \kappa^0 + \kappa^1 + \kappa^2$ with 
\begin{eqnarray*}
\begin{split}
\kappa^0(x_1,x_2; p_0) & = 8+p_0,\\
\kappa^1(x_1,x_2; p_1) & = e^{x_1+x_2+p_1}\cos \left (\frac{2\pi x_2}{\varepsilon} \right ) \sin \left (\frac{2\pi x_1}{\varepsilon} \right ), \\
\kappa^2(x_1,x_2; p_2) & = e^{x_1 x_2+p_2}\cos \left (\frac{2\pi x_1}{\varepsilon} \right ) \sin \left (\frac{2\pi x_2}{\varepsilon} \right ).
\end{split}
\end{eqnarray*}
Here, we set $\varepsilon = 0.1$. The parameters $p_0$, $p_1$, and $p_2$ are of uniform distribution with 
$$ p_0 \sim U[-2,2], \quad p_1 \sim U[-1.2, 1.2], \quad \text{and} \quad p_2 \sim U[-1.5, 1.5].$$

There are three random parameters associated with the problem. However, we will design a two-stage training process 
since we have two different scales only. Our experiments show that an addition third features will merely improve the result of training. The training process can be formulated as follows:
\begin{itemize}
\item In the first stage, we solve the following optimization problem
$$
\min_{G_1} \norm{ (u_{H, i}) - G_1(D_1(\kappa(x, s) ) }_1,
$$
where $D_1(\kappa(x, s)) = (f_i^0) \in\mathbb{R}^N$ and $G_1$ is a generator in the first stage. Here, one can set  $G_1$ as a two-layer fully connected network.
\item In the second stage, we minimize 
$$
\min_{G_2, C_2} \norm{ (u_{H, i}) - C_2(u_{H, s_{1}}^c, u_{H, s_{1}} )  }_1,
$$
where $u_{H, s_{1}} = G_2(D_2( \kappa(x, s) )) \in\mathbb{R}^N$ and $u_{H, s_{1}}^c \in\mathbb{R}^N$ is the prediction generated from the previous stage. We remark that $D_2(\kappa(x, s)) = (f_i^1)$ and $G_2$ has identical structure as $D_1$. One can set $C_2: \mathbb{R}^l\times\mathbb{R}^l\rightarrow \mathbb{R}^l$ to have two fully connected layers, which combines and corrects the predictions.
\end{itemize}

\subsection{Numerical Results}
In this section, we present some numerical results using the proposed multi-stage learning strategy for the steady-state model. 

We generate $1200$ training samples and test on $400$ samples. The relative errors of the prediction using the mean of the training data is about $0.33$.
We first train the network with $(f_i^0)$ and $(f_i^1)$ in the first and second stages, respectively. 
The results of the 2-stage training is shown in Table \ref{kappa_example}. One can clearly observe from the results that the second stage helps enhance the accuracy of the approximation generated from the previous stage. 

\begin{table}[H]
\centering
\begin{tabular}{||c c||} 
\hline
Stage 1 Error & Stage 2 Error \\ [0.5ex] 
\hline
0.22485 & 0.10643  \\ [0.5ex]
\hline
\end{tabular}
\caption{$L^2$ errors of each stage for steady-state model.}
\label{kappa_example}
\end{table}

On the other hand, one can also apply the idea of pooling which is in Section \ref{exp3} to upgrade the two-stage training. The guaranteed benefit of pooling is the save on computing the multiscale features. The result is shown in Table \ref{kappa_example_pool}. The result is very closed to the previous experiment which we use multiscale features in both stages and we claim that this is a good update.

\begin{table}[H]
\centering
\begin{tabular}{||c c||} 
\hline
Stage 1 Error & Stage 2 Error \\ [0.5ex] 
\hline
0.19730 & 0.10842  \\ [0.5ex]
\hline
\end{tabular}
\caption{$L^2$ errors of each stage for steady-state model with max-pooling in Stage 2.}
\label{kappa_example_pool}
\end{table}

\section{Conclusion} \label{conclusion}
In this research, we have proposed a multi-stage training methodology coupling with multiscale model reduction for multiscale problems. 
In the first stage of the training process, a rough prediction to the model is generated. Then, in the following stages, one can iteratively improve the accuracy of the prediction. Each stage shares almost identical network structure. 
We have verified our method on two time-dependent linear and nonlinear PDEs. With brief modification, we have numerically shown that our method is applicable for steady-state model problem. Several upgrades of the methodology have been also discussed.

In the future, we will study the design of a good reduced order model and how to recovered a fine model from the learnt reduced order model.  The layers information of the reduced order model can also be studied in the light of the improvement of each stage.

\section*{Acknowledgement}

The research of Eric Chung is partially supported by the Hong Kong RGC General Research Fund (Project numbers 14304719 and 14302018) and CUHK Faculty of Science Direct Grant 2019-20.

\bibliographystyle{plain}

\bibliography{references}

\end{document}